\documentclass[12pt]{article}
\usepackage{amsfonts}

\usepackage{latexsym}
\usepackage{amssymb}
\usepackage{epsfig}
\usepackage{amsmath}
\usepackage{amsthm}
\usepackage[mathscr]{eucal}
\usepackage{subfigure}
\usepackage{graphicx}

\newtheorem{Theorem}{Theorem}
\newtheorem{Proposition}{Proposition}

\renewcommand{\qed}{\hfill{\ \ \rule{2mm}{2mm}} \vspace{0.2in}}

\newcommand{\ind}{1\hspace{-2.3mm}{1}}

\begin{document}

\title{Rates of linear codes with low decoding error probability}
\author{ \textbf{Ghurumuruhan Ganesan}
\thanks{E-Mail: \texttt{gganesan82@gmail.com} } \\
%EndAName
\ \\
New York University, Abu Dhabi}
\date{}
\maketitle

\begin{abstract}
Consider binary linear codes obtained from bipartite graphs as follows. There are~\(k \geq 1\) left nodes each representing a message bit and there are~\(m = m(k)\) right nodes each representing a parity bit, generated from the corresponding set of message node neighbours. Both the message and the parity bits are sent through a memoryless binary input channel that either retains, flips or erases each transmitted bit, independently. Based on the received set of symbols, the decoder at the receiver obtains an estimate of the original message sent. If the decoding error probability~\(P_k \longrightarrow 0\) and the average degree per parity node remains bounded as~\(k \rightarrow \infty,\) then the rate of the code~\(\frac{k}{k+m} \longrightarrow 0\) as~\(k \rightarrow \infty.\) %The analysis and the results also extend to binary erasure channel with positive erasure probability.

\vspace{0.1in} \noindent \textbf{Key words:} Linear codes, low decoding error probability, asymptotic rates.

\vspace{0.1in} \noindent \textbf{AMS 2000 Subject Classification:} Primary:
60J10, 60K35; Secondary: 60C05, 62E10, 90B15, 91D30.
\end{abstract}

\bigskip

\setcounter{equation}{0}
\renewcommand\theequation{\thesection.\arabic{equation}}
\section{Introduction} \label{intro}
Parity check codes are used extensively in today's communication systems particularly in the form of Low Density Parity Check (LDPC) Codes (see~\cite{amin} for an introduction). One of the main challenges here is to achieve low decoding error probability. Previous papers have mainly focused on decoding schemes that achieve low error probability (see for example~\cite{luby} and references therein). The emphasis there is to design schemes that achieve low error probability but possibly at the cost of increased overhead.

In this paper, we study the rate versus decoding error probability tradeoff and show that low decoding error probability necessarily requires a low rate or equivalently a large number of parity bits to be appended to the message. In other words, if the decoder is such that the asymptotic decoding error probability converges to zero as the number of message bits~\(k \rightarrow \infty,\) then the asymptotic encoded rate also converges to zero as~\(k \rightarrow \infty.\)

%Let~\(X = (X_1,\ldots,X_k)\) be the random binary message bits, independent and identically distributed (i.i.d) with~
%For~\(1 \leq j \leq m,\) let~\(R_k(j)\) be the set of message nodes adjacent to the parity node~\(j\) an define
%to be the~\(j^{th}\) parity bit.

\subsection*{Model description}
We are interested in sending a random message through a communication channel reliably. We describe the underlying communication system below.

\subsubsection*{Messages}
Messages are~\(k-\)bit vectors satisfying the following condition:\\
\((A1)\) A random message~\(X = (X_1,\ldots,X_k)\) has independent and identically distributed (i.i.d.) bits~\(X_i \in \{0,1\}\) with
\begin{equation}\label{x_val}
\mathbb{P}(X_i = 0) = p_x = 1-\mathbb{P}(X_i = 1).
\end{equation}
In particular, this implies that the \emph{raw rate} defined as
\begin{equation}\label{r_raw}
R_{raw} := \frac{H(X)}{k} = H(p_x) > 0,
\end{equation}
where~\(H(p_x) = -p_x\log{p_x} - (1-p_x)\log(1-p_x)\) and
\begin{equation}\label{entro}
H(X) := -\sum_{w} p(w)\log{p(w)}
\end{equation}
is the entropy of the vector~\(X\) (see Chapter~\(1,\) Section~\(1.1\) of~\cite{joy}). In~(\ref{entro}),~\(p(.)\) is the probability mass function of~\(X\) and the summation is over all possible~\(k-\)bit vectors. All logarithms are to the base~\(2\) and for simplicity we assume throughout that~\(p_x = \frac{1}{2}\) so that~\(H(p_x) = 1.\)

\subsubsection*{Encoder}
We consider binary linear codes obtained from bipartite graphs as follows. There are~\(k \geq 1\) left nodes called \emph{message nodes} and there are~\(m = m(k)\) right nodes called \emph{parity nodes}. For parity node~\(1 \leq j \leq m,\) let~\(R_k(j)\) be the message nodes adjacent to~\(j.\) The~\(j^{th}\) parity bit \(Z_j\) is obtained as
\begin{equation}\label{z_i_def}
Z_j = \oplus_{w \in R_k(j)} X_w,
\end{equation}
where~\(\oplus\) is XOR operation, i.e., addition modulo~\(2.\) The vector~\[(X,Z_1,\ldots,Z_m) = (X_1,\ldots,X_k,Z_1,\ldots,Z_m)\] is the \emph{codeword} associated with the message~\(X\) and the \emph{encoded rate} is defined as
\begin{equation}\label{r_enc}
R_{enc} := \frac{H(X)}{k+m} = \frac{k}{k+m},
\end{equation}
by~(\ref{r_raw}). We make the following assumption regarding the encoder:\\
\((A2)\) For~\(1 \leq j \leq m,\) let~\(\#R_k(j)\) be the degree of the parity node~\(j\) and suppose that the average degree per parity node remains bounded as~\(k \rightarrow \infty;\) i.e.,
\begin{equation}\label{f_bar_bd}
\limsup_k \frac{1}{m}\sum_{j=1}^{m}\#R_k(j) < \infty.
\end{equation}

\subsubsection*{Channel}
The codeword~\((X,Z_1,\ldots,Z_m)\) is sent through a binary input channel which introduces noise that either retains, flips or erases the transmitted bit.
Formally, we assume that the noise alphabet is~\(\{\alpha_0,\alpha_1,\alpha_{er}\}\) and the~\(i^{th}\) received message symbol is
\begin{equation}\label{n_def_bec}
\tilde{X}_i = \ind(N_x(i) = \alpha_{er})\alpha_{er} + \ind(N_x(i) = \alpha_0) X_i + \ind(N_x(i) = \alpha_1)(1-X_i)
\end{equation}
Here~\(\alpha_{er}\) is the erasure symbol and~\(\tilde{N}_x(i)\) is the noise symbol. Similarly, the~\(j^{th}\) received parity symbol is
\begin{equation}\label{n_def_bec2}
\tilde{Z}_j = \alpha_{er}\ind(N_z(j) = \alpha_{er}) + \ind(N_z(j) = \alpha_0) Z_j + \ind(N_z(j) = \alpha_1)(1-Z_j).
\end{equation}
The overall received codeword is
\begin{equation}\label{y_def}
Y = (\tilde{X}_1,\ldots,\tilde{X}_k,\tilde{Z}_1,\ldots,\tilde{Z}_m).
\end{equation}
\((A3)\) We assume that the noise random variables~\(\{N_x(i)\}\) and~\(\{N_z(j)\}\) are independent and identically distributed (i.i.d.) with
\begin{equation}\label{noise_bec}
\mathbb{P}(N_x(i) = \alpha_{er}) = p_{er}, \mathbb{P}(N_x(i) = \alpha_1) = p_1 \text{ and }\mathbb{P}(N_x(i) = \alpha_0) =1-p_1-p_{er}.
\end{equation}
The term~\(p_{er} + p_1\) is the probability that a channel error occurs; i.e., the noise in the channel corrupts (either erases or flips) a transmitted bit. We also assume that the noise is independent of the transmitted bits~\(\{X_i\}\) and~\(\{Z_j\}.\)

%The above channel model is general and with particular choices of~\(\alpha_0,\alpha_1\) and~\(\alpha_{er},\) we realize various channels. For example if~\(\alpha_0 =0,\alpha_1 = 1\) and~\(p_0 = p = 1-p_1,\) then we obtain the binary symmetric channel (BSC). Similarly if~\(p_0 =p = 1-p_{er},\) we obtain the binary erasure channel (BEC). If we let~\(p_{1} = p\)

\subsubsection*{Decoder}
At the receiver, a pre installed decoder uses the received word~\(Y\) to obtain an estimate~\(\hat{X}\) of the message sent and let
\begin{equation}\label{dec_error}
P_k = \mathbb{P}(X \neq \hat{X})
\end{equation}
be the decoding error probability. The following is the main result of this paper.
\begin{Theorem}\label{thm1} Suppose assumptions~\((A1)-(A3)\) hold.
If the decoding error probability~\(P_k \longrightarrow 0\) as~\(k \rightarrow \infty,\) then the encoded rate~\(R_{enc} = \frac{k}{k+m} \longrightarrow 0\) as~\(k \rightarrow \infty.\)
\end{Theorem}
In other words, any code having low decoding error probability must necessarily contain a lot of parity bits. One example of such a code is the~\(r-\)repetition code, where each message bit is simply repeated~\(r\) times. Recall that for a fixed~\(r,\) an~\(r-\)repetition code has an encoded rate of~\(\frac{1}{r+1}\) and using majority decision rule, it is possible to correct up to~\(\frac{r-1}{2}\) channel errors, irrespective of the number of bits~\(k\) in the message (see~\cite{wiki}).  If however, we allow~\(r = r(k)\) to depend on~\(k,\) we can correct all errors in the message with high probability. %CHECK IF WE NEED TO PUT IT IN OR OUT...PRKVMM+etC..

\begin{Proposition}\label{prop1} Suppose~\(2p_1+p_{er} < 1\) and \(r = r(k) = M\log{k}.\) There are constants~\(M_0 = M_0(p_1,p_{er}) \geq 1\) and~\(K_0 = K_0(p_1,p_{er}) \geq 1\) so that the following holds for all~\(M \geq M_0\) and~\(k \geq K_0:\) For an~\(r-\)repetition code, the decoding error probability with the majority decision rule is bounded above by~\(P_k \leq \frac{1}{k}.\)
\end{Proposition}

%see prkvmm+etc... WRT MORE +etC...

The paper is organized as follows. In Section~\ref{pf1}, we prove Theorem~\ref{thm1} and Proposition~\ref{prop1}. %for BSC and in Section~\ref{pf2}, we prove Theorem~\ref{thm1} for BEC.

\section*{Proof of Theorem~\ref{thm1} and Proposition~\ref{prop1}}\label{pf1}
Recall that~\(X\) is the message and~\(Y\) as defined in~(\ref{y_def}) is the received codeword. Define
\begin{equation}
H(X|Y) := -\sum p(x,y)\log{p(x|y)}
\end{equation}
to be the uncertainty in~\(X\) given the random vector~\(Y,\) where~\(p(x,y)\) and~\(p(x|y)\) respectively, refer to probability mass functions of the joint distribution of~\((X,Y)\) and the conditional distribution of~\(X\) given~\(Y\) (see Chapter~\(1,\)~\cite{joy}). Since the total number of messages is~\(2^{k},\) we have from Fano's inequality (Theorem~\(2.10.1,\)~\cite{joy}) that
\begin{equation}\nonumber
H(X|Y) \leq H(X|\hat{X}) \leq H(P_k) + P_k\log\left(2^{k} - 1\right) \leq 1 + kP_k
\end{equation}
and so
\begin{equation}\label{eq1}
\frac{1}{k}H(X|Y) \leq \frac{1}{k} + P_k \longrightarrow 0
\end{equation}
as~\(k \rightarrow \infty.\)

%H(X_1|Y) + H(X_2|Y,X_1) + \ldots + H(X_n|Y,X_1,\ldots,X_{n}).

To evaluate~\(H(X|Y),\) let~\(X_0 = 0\) and write
\begin{equation}
H(X|Y)  = \sum_{i=1}^{k}  H(X_i|Y,X_1,\ldots,X_{i-1}) \geq \sum_{i=1}^{k}H(X_i|\tilde{X}_i,\{\tilde{Z}_j\},\{X_w\}_{w \neq i}). \label{succ}
\end{equation}
The first equality in~(\ref{succ}) follows by chain rule for entropy (Theorem~\(2.5.1\),~\cite{joy}) and the inequality in~(\ref{succ}) follows from the data processing inequality (Theorem~\(2.8.1,\)~\cite{joy}).

We evaluate each term in the summation in~(\ref{succ}) separately. First, we use the received parity symbols~\(\tilde{Z}_1,\ldots,\tilde{Z}_m\) to obtain estimates for the~\(i^{th}\) transmitted bit~\(X_i.\) Formally, for~\(1 \leq i \leq k\) let~\(T_k(i)\) denote the set of parity nodes adjacent to the message node~\(i.\) Recall that for~\(u \in T_k(i),\) the term~\(R_k(u)\) denotes the set of message nodes adjacent to the parity node~\(u\) and by definition~\(i \in R_k(u).\) For~\(1 \leq i \leq k,\) define
\begin{eqnarray}
\hat{X}_i(u) &:=& \alpha_{er} \ind(\tilde{Z}_u = \alpha_{er}) + \ind(\tilde{Z}_u \neq \alpha_{er}) \tilde{Z}_u \oplus_{w \in R_k(u) \setminus \{i\}}X_w \nonumber\\
&=& \alpha_{er} \ind(N_z(u) = \alpha_{er}) + \ind(N_z(u) = \alpha_0)X_i + \ind(N_z(u) = \alpha_1) (1-X_i). \nonumber\\\label{hat_def_bec}
\end{eqnarray}
Equation~(\ref{hat_def_bec}) follows from the expression for~\(\tilde{Z}_u\) in~(\ref{n_def_bec2}) and the fact that if~\(Z_u = X_i \oplus_{w \in R_k(u) \setminus \{i\}}X_w,\) then~\[1-Z_u = (1-X_i) \oplus_{w \in R_k(u) \setminus \{i\}}X_w.\] The map~\[\left(\tilde{X}_i,\{\tilde{Z}_j\},\{X_w\}_{w \neq i}\right) :\longrightarrow \left(\tilde{X}_i,\{\hat{X}_i(u)\}_{u \in T_k(i)}, \{\tilde{Z}_j\}_{j \notin T_k(i)},\{X_w\}_{w \neq i}\right)\] is one to one and invertible and so the \(i^{th}\) term in the final summation in~(\ref{succ}) is
\begin{equation}\label{map4}
H(X_i|\tilde{X}_i,\{\tilde{Z}_j\},\{X_w\}_{w \neq i}) = H\left(X_i|\tilde{X}_i,\{\hat{X}_i(u)\}_{u \in T_k(i)}, \{\tilde{Z}_j\}_{j \notin T_k(i)},\{X_w\}_{w \neq i}\right).
\end{equation}

The set of random variables~\((\{\tilde{Z}_j\}_{j \notin T_k(i)},\{X_w\}_{w \neq i})\) are independent of the rest of random variables~\((\tilde{X}_i,\{\hat{X}_i(u)\}_{u \in T_k(i)})\) and are also independent of~\(X_i.\) Thus
\begin{equation}\label{map5}
H\left(X_i|\tilde{X}_i,\{\hat{X}_i(u)\}_{u \in T_k(i)}, \{\tilde{Z}_j\}_{j \notin T_k(i)},\{X_w\}_{w \neq i}\right) = H\left(X_i|\tilde{X}_i,\{\hat{X}_i(u)\}_{u \in T_k(i)}\right)
\end{equation}
and substituting this into~(\ref{succ}) gives
\begin{equation}\label{map6}
H(X|Y) \geq \sum_{i=1}^{k} G(d_k(i))
\end{equation}
where~\(d_k(i) := \#T_k(i)\) is the degree of the message node~\(i\) and~\[G(d_k(i)) := H\left(X_i|\tilde{X}_i,\{\hat{X}_i(u)\}_{u \in T_k(i)}\right) > 0\] is the uncertainty in the bit~\(X_i\) given~\(d_k(i)+1\) independently noise corrupted copies.

We have the following properties regarding~\(G(.).\)\\
\((g1)\) Using the fact that conditioning reduces entropy, we obtain that~\(G(d)\) is a decreasing function of~\(d.\)\\
\((g2)\) Using~(\ref{map6}) and~(\ref{eq1}) we get that
\begin{equation}\label{map7}
\frac{1}{k}\sum_{i=1}^{k} G(d_k(i)) \longrightarrow 0
\end{equation}
as~\(k \rightarrow \infty.\)

We use properties~\((g1)-(g2)\) to get the following properties.\\
\((g3)\) The average degree per message node
\begin{equation}\label{eq_ave_left}
\frac{1}{k}\sum_{i=1}^{k}d_k(i) \longrightarrow \infty
\end{equation}
as~\(k \rightarrow \infty.\) \\
\((g4)\) The encoded rate~\(\frac{k}{k+m} \longrightarrow 0\) as~\(k \rightarrow \infty.\)\\
This proves Theorem~\ref{thm1}.

\emph{Proof of~\((g3)-(g4)\)}: We prove~\((g3)\) first. For integer~\(q \geq 1,\) let
\begin{equation}\label{s_def}
S_k(q) = \{i :  d_k(i) \leq q\}
\end{equation}
be the set of message nodes whose degree is at most~\(q.\) For a fixed~\(q,\) it is true that
\begin{equation}\label{s_conv}
\frac{\#S_k(q)}{k} \longrightarrow 0
\end{equation}
as~\(k \rightarrow \infty.\) If~(\ref{s_conv}) is not true, then there exists~\(\epsilon_0 > 0\) and a subsequence~\(\{k_r\}\) such that~\(\frac{\#S_{k_r}(q)}{k_r} \geq \epsilon_0\) for all large~\(r.\) Using property~\((g1)\) that~\(G(.)\) is decreasing, we get that
\[\frac{1}{k_r}\sum_{i=1}^{k_r} G(d_{k_r}(i)) \geq \frac{1}{k_r}\sum_{i \in S_{k_r}(q)} G(d_{k_r}(i)) \geq G(q)\frac{\#S_{k_r}(q)}{k_r} \geq \epsilon_0 G(q) >0\] for all large~\(r,\) contradicting~(\ref{map7}) in property~\((g2).\)

From the above paragraph, we obtain that~(\ref{s_conv}) is true and so for any integer~\(q \geq 1,\) we get that
\[\frac{1}{k}\sum_{i=1}^{k}d_k(i) \geq \frac{1}{k}\sum_{i \notin S_k(q)} d_k(i) \geq q\left(\frac{k-S_k(q)}{k}\right) \geq \frac{q}{2}\] for all large~\(k.\) Since~\(q \geq 1\) is arbitrary, we get~(\ref{eq_ave_left}).

To prove~\((g4),\) we use the fact that the number of edges in the graph is~\[\sum_{i=1}^{k}d_k(i) = \sum_{j=1}^{m} f_k(j)\] where~\(f_k(j) = \#R_k(j)\) is the degree of the parity node~\(j.\) Using~\((g3),\) we therefore get~
\begin{equation}\label{eq_left_right}
\frac{1}{k}\sum_{j=1}^{m} f_k(j) = \frac{m}{k}\frac{1}{m}\sum_{j=1}^{m}f_k(j) \longrightarrow \infty
\end{equation}
as~\(k \rightarrow \infty.\) Since by assumption, the average degree per parity node is bounded~(see (\ref{f_bar_bd})) we get from~(\ref{eq_left_right}) that~\(\frac{m}{k} \longrightarrow \infty\) and so~\(\frac{k}{k+m} \longrightarrow 0\) as~\(k \rightarrow \infty.\)~\(\qed\)

\emph{Proof of Proposition~\ref{prop1}}: Let~\(X = (X_1,\ldots,X_k)\) be the message bits. For~\(1 \leq i \leq k\) and \(1 \leq j \leq r,\) define~\(Z_i(j) = X_i\) be the parity bits for the message bit~\(X_i.\) Thus each message bit is repeated~\(r\) times and for convenience define~\(Z_i(0) = X_i\) to be the message bit to be transmitted. Let~\(\{\tilde{Z}_j(i)\}\) be corresponding received symbols as defined in~(\ref{n_def_bec2}).

The decoding is majority based as follows. For each~\(1 \leq i \leq k\) and~\(l \in \{0,1\},\) let~\(W_l(i) \subseteq \{0,1,2,\ldots,r\}\) be the random set of all indices for which the received symbol is~\(l;\) i.e.,~\[\tilde{Z}_i(j) =0 \text{ for all }j \in W_0(i)\text{ and }\tilde{Z}_i(j) = 1 \text{ for all } j \in W_1(i).\] If~\(\#W_1(i) \geq \#W_0(i),\) set~\(\hat{X}_i = 1;\) else set~\(\hat{X}_i = 0.\) The estimated message is~\(\hat{X} = (\hat{X}_1,\ldots,\hat{X}_k).\)

A decoding error occurs if~\(\hat{X}_i \neq X_i\) for some~\(1 \leq i \leq r.\) For a fixed~\(1 \leq i \leq k\) and~\(0 \leq j \leq r,\) let~\(N_z(i,j) \in \{0,1,\alpha\}\) be the noise random variable affecting the bit~\(Z_i(j)\) as in~(\ref{n_def_bec2}). Message bit~\(i\) is decoded wrongly if and only if~\[\sum_{j=0}^{r}\ind(N_z(i,j)=1) \geq \sum_{j=0}^{r}\ind(N_z(i,j) = 0).\] Defining~\[L(i,j) = \ind(N_z(i,j)=1) -\ind(N_z(i,j) = 0) \in \{-1,1\}\]  we have that~\[\mathbb{E}L(i,j) = p_1-(1-p_1-p_{er})  = 2p_1 + p_{er} - 1 <0,\] by the assumption in the statement of the Proposition.

For a fixed~\(1 \leq i \leq k,\) the random variables~\(\{L(i,j)\}_{0 \leq j \leq r}\) are i.i.d and so using the Chernoff bound, we have for~\(s  >0\) and~\(c \geq 0\) that
\begin{equation}\label{lij_est}
\mathbb{P}\left(\sum_{j=0}^{r} L(i,j) \geq c\right) \leq e^{-sc}\prod_{j=0}^{r}\mathbb{E}e^{sL(i,j)} = e^{-sc} \left(e^{s}p_1 + e^{-s}(1-p_1-p_\alpha)\right)^{r+1}.
\end{equation}
Writing~\(e^{s} = 1+s + R_1(s)\) and~\(e^{-s} = 1-s + R_2(s),\) we have
\[p_1e^{s} + e^{-s}(1-p_1-p_{er}) = 1-(1-2p_1-p_{er})s + T(s),\] where~\(T(s) = R_1(s)p_1 + R_2(s)(1-p_1-p_{er}).\) Choosing~\(s > 0\) small, we have~\(|T(s)| \leq s^2\) and~\(1-(1-2p_1-p_{er})s + T(s)\leq \delta\) for some constant~\(\delta < 1.\) Substituting into~(\ref{lij_est}) and setting~\(c = 0\) gives~\[\mathbb{P}\left(\sum_j L(i,j) \geq 0\right) \leq \delta^{r+1} \leq \frac{1}{k^2}\] if~\(r = \frac{2}{\delta}\log{k}.\) But~\(\sum_j L(i,j) \geq 0\) if and only if the bit~\(X_i\) is decoded wrongly i.e.,~\(\hat{X}_i \neq X_i\) and so~\(\mathbb{P}(\hat{X}_i \neq X_i) \leq \frac{1}{k^2}\) and so the overall decoding error probability is at most~\(\frac{1}{k} \longrightarrow 0\) as~\(k \rightarrow \infty.\)~\(\qed\)

%\subsection*{Acknowledgements}
%I thank Professor Federico Camia for crucial comments and for my fellowships.

%WRIT MORE +ETC.. FOR BEC ALSO AND ALSO POLISH +eTC...PARKVMM +eTC FOR OUR BNFT +eTCC....

%erere

\setcounter{equation}{0} \setcounter{Lemma}{0} \renewcommand{\theLemma}{II.%
\arabic{Lemma}} \renewcommand{\theequation}{II.\arabic{equation}} %
\setlength{\parindent}{0pt}

%\section*{Appendix~II: Proof of (\ref{p_t_y_ineq1}), (\ref{e_2_m_c1}) and (\ref{e_2_m_c})}\label{app1}
%\begin{figure}
%\centering
%\includegraphics[width=3.5in]{cap_scheme.eps}
%\caption{Modified coding scheme of~\cite{Goldsmith}.} \label{fig_eta_res}
%\end{figure}

%\begin{figure}
%\centering
%\includegraphics[width=4in]{cap_vs_se.eps}
%\caption{Capacity as a function of network population \(n\) for various values of \(\sigma_e.\)} \label{cap_se}
%\end{figure}

%\bibliographystyle{IEEEtran}
%\bibliography{IEEEabrv,ref_mob}

% biography section
%
% If you have an EPS/PDF photo (graphicx package needed) extra braces are
% needed around the contents of the optional argument to biography to prevent
% the LaTeX parser from getting confused when it sees the complicated
% \includegraphics command within an optional argument. (You could create
% your own custom macro containing the \includegraphics command to make things
% simpler here.)
%\begin{biography}[{\includegraphics[width=1in,height=1.25in,clip,keepaspectratio]{mshell}}]{Michael Shell}
% where an .eps filename suffix will be assumed under latex, and a .pdf suffix
% will be assumed for pdflatex; or if you just want to reserve a space for
% a photo:

% You can push biographies down or up by placing
% a \vfill before or after them. The appropriate
% use of \vfill depends on what kind of text is
% on the last page and whether or not the columns
% are being equalized.

%\vfill

% Can be used to pull up biographies so that the bottom of the last one
% is flush with the other column.
%\enlargethispage{-5in}

% that's all folks

\bibliographystyle{plain}

\end{document}